\documentclass[11pt]{article}

\pdfoutput=1

\usepackage{latexsym}
\usepackage{amssymb}
\usepackage{amsfonts}
\usepackage{graphicx}

\def\T{{ \mathrm{\scriptscriptstyle T} }}

\newcommand{\bSigma}{\mbox{\boldmath$\Sigma$}}

\newcommand{\bA}{{\mathbf A}}
\newcommand{\bI}{{\mathbf I}}
\newcommand{\bM}{{\mathbf M}}
\newcommand{\bP}{{\mathbf P}}
\newcommand{\bR}{{\mathbf R}}
\newcommand{\bS}{{\mathbf S}}
\newcommand{\bT}{{\mathbf T}}
\newcommand{\bU}{{\mathbf U}}
\newcommand{\bV}{{\mathbf V}}
\newcommand{\bW}{{\mathbf W}}
\newcommand{\bX}{{\mathbf X}}
\newcommand{\bt}{{\mathbf t}}
\newcommand{\bu}{{\mathbf u}}
\newcommand{\bv}{{\mathbf v}}
\newcommand{\bzero}{{\mathbf 0}}

\begin{document}

\vskip 1 cm

\begin{center}
{\LARGE
Statistics of ambiguous rotations
}
\end{center}

\medskip

\begin{center}
R. Arnold, \\
{\small  School of Mathematics and Statistics,
Victoria University of Wellington,  \\
PO Box 600, Wellington, New Zealand,}  \\
 \smallskip
 P.E. Jupp,  \\
{\small School of Mathematics and Statistics, University of St Andrews,  \\
St Andrews, Fife KY16 9SS, UK,}  \\
 \smallskip
   H. Schaeben, \\
{\small Geophysics and Geoscience Informatics,
TU Bergakademie Freiberg,  \\
Germany}
\end{center}
    
\medskip

\vskip 1cm

\begin{abstract}

The orientation of a rigid object can be described by a rotation 
that transforms it into a standard position.
For a symmetrical object the rotation is known only up to multiplication by an element of the 
symmetry group.
Such ambiguous rotations arise in biomechanics, crystallography and seismology. We develop methods for 
analyzing data of this form.
A test of uniformity is given. 
Parametric models for ambiguous rotations are presented, tests of location are 
considered, and a regression model is proposed.  
A brief illustrative example involving orientations of diopside crystals is given.
\end{abstract}


\section{Introduction}

Data that are rotations of $\mathbb{R}^3$ occur in various areas of science, such as 
 palaeo-magnetism (Pesonen et al., 2003; Villala\'{\i}n et al., 2016;  Koymans et al., 2016), 
 plate tectonics and
 seismology (Stein \& Wysession, 2003; Hardebeck, 2006; Arnold \& Townend, 2007;
 Walsh et al., 2009; Khalil \& McClay, 2016), 
 biomechanics (Rivest, 2005; Lekadir et al., 2015; Spronck et al., 2016), 
 crystallography (Hielscher et al., 2010; Griffiths et al., 2016)
 and texture analysis, i.e., analysis of orientations of crystalites
(Kunze \& Schaeben, 2004, 2005; Du et al., 2016).  
The sample space is the 3-dimensional rotation group, $SO(3)$,
and methods for handling such data are now an established part of directional statistics;
see \S13.2 of Mardia \& Jupp (2000).
In some contexts the presence of symmetry means that 
the rotations are observed subject to ambiguity, 
so that it is not possible to distinguish a rotation $\bX$ from $\bX \bR$ for any rotation $\bR$ in some 
given subgroup $K$ of $SO(3)$.
From the mathematical point of view,  the sample space is the quotient $SO(3)/K$ of $SO(3)$
 by $K$.
 Such spaces arise in many scientific contexts:
the case in which $K$ is generated by the rotations through $180^{\circ}$ about the coordinate axes 
gives the orthogonal axial frames considered by   
Arnold \& Jupp (2013), which can be used to describe aspects of earthquakes;
many groups $K$ of low order occur as the symmetry groups of crystals;     
the icosahedral group is the symmetry group of 
some carborane molecules (Jemmis, 1982),  of most closed-shell viruses (Harrison, 2013),
of the natural quasicrystal, icosahedrite  (Bindi et al., 2011),
and of the blue phases of some liquid crystals 
(Seideman,1990, \S 6.1.2).      
The object of this paper is to give a unified account of 
some general tools for the analysis of data consisting of ambiguous rotations with 
a finite symmetry group.

\section{Ambiguous rotations}

\subsection{Symmetry groups}

The orientation of a rigid object in $\mathbb{R}^3$ can be described by a rotation 
that transforms it into some standard position.
If the object is asymmetrical then this rotation is unique, so that the orientations of the object 
correspond to elements of the rotation group $SO(3)$.
If the object is symmetrical then the set of rotations that have no visible effect on the object
forms a subgroup $K$ of $SO(3)$.
Then the orientations of the object 
correspond to elements of the homogeneous space $SO(3)/K$,
i.e.\ the set of equivalence classes of elements of $SO(3)$ under the right action of $K$.
We shall consider the cases in which $K$ is finite.
In particular, the orientations of T-shaped, X-shaped and $+$-shaped objects in $\mathbb{R}^3$ are elements of $SO(3)/K$ with 
$K = C_2, D_2$ and $D_4$, respectively.
For $\bU$ in $SO(3)$ we shall denote the equivalence class of $\bU$ in
$SO(3)/K$ by $[\bU]$.

The finite subgroups of $SO(3)$ are known also as the point groups of the first kind. The classification result for these groups, given e.g.\ in Miller (1972),    
states that any such group is isomorphic to one of the
following: the cyclic groups, $C_r$, for $r = 1, 2, \dots$, 
the dihedral groups, $D_r$, for $r = 2, 3, \dots$, 
the tetrahedral group, $T$, the octahedral group, $O$,
and the icosahedral group, $Y$. 
These groups are listed in Table 1,   together with the frames of 
vectors that will be used to represent elements of the sample spaces $SO(3)/K$. 
The group $C_1$ has one element, the identity, $\bI_3$.

\subsection{Frames and symmetric frames}

For each point group $K$ of the first kind,
every element of $SO(3)/K$ can be represented uniquely by a $K$-frame, i.e.\ 
an equivalence class of a frame, meaning a set of vectors or axes in $\mathbb{R}^3$.
For $K = C_r$  with $r \ge 3$ or $K = D_r$ with $r \ge 3$, it is convenient to take the 
vectors of the frame to be unit normals to the sides of a regular $r$-gon;
for $K = C_2$ we take a unit vector and an axis orthogonal to it; 
for $K = D_2$ we take a pair of orthogonal axes; 
for $K = T$, $O$ or $Y$, it is convenient to take the 
vectors to be unit normals to the sides of a regular tetrahedron, cube or dodecahedron, 
respectively.
Permutation of the vectors of the frame by the action of $K$ leads to ambiguity. 
This ambiguity is removed by passing to the corresponding $K$-frame, i.e.\ the equivalence class of the frame under such permutations.
The $K$-frames will be denoted by square brackets, e.g.\ for $K = C_r$,
$[\bu_1,\dots ,  \bu_r ]$ denotes the $K$-frame arising from  $(\bu_1,\dots ,  \bu_r)$.
By a symmetric frame, we shall mean a $K$-frame for some $K$.
The frames that we consider are listed in Table  1, together with an indication of the ambiguities.

\begin{table}[ht]
\begin{center}
\caption{Symmetry groups and frames.} 
{\scriptsize  
\begin{tabular}{llll}
\hline 
	Group     & Name  &  Frame  &   \\
\hline 
  $C_1$       &   trivial   &    $( \bu_ 1 , \bu_2, \bu_3)$      &  
$\bu_1,  \bu_2, \bu_3$ orthonormal, $\bu_ 3 = \bu_1 \times \bu_2$    \\  
    $C_2$       &   cyclic   &      $( \bu_ 0 , \pm \bu_1)$      &  
$\bu_0, \bu_1$ orthonormal    \\  
    $C_r \quad  (r \ge 3)$  &   cyclic           & $( \bu_1 ,  \dots , \bu_r )$  &  
$\bu_1 ,  \dots , \bu_r$ coplanar,  \\ 
              &            &              & known up  to cyclic order,   \\  
              &            &              & $\bu_i ^{\T} \bu_{i-1} = \cos  (2 \pi /r)$ for $i = 2, \dots, r$   \\  
   $D_2$       &   dihedral   &      $( \pm \bu_1 , \pm \bu_2)$      &  orthogonal axes    \\  
$D_r \quad  (r \ge 3)$   & dihedral  & $( 
\bu_1 ,  \dots , \bu_r )$    & $\bu_1 ,  \dots , \bu_r$ coplanar,  \\ 
              &            &             & known up  to cyclic order and reversal,   \\  
            &            &                 &$\bu_i ^{\T} \bu_{i-1} = \cos  (2 \pi /r)$ for $i = 2, \dots, r$   \\    
 $T = A_4$   & tetrahedral   & 
   $\{   \bu_1,\dots , \bu_4  \}$  &  $\bu_i ^{\T} \bu_j = -1/3$ for $ i \ne j$   \\    
$O = \Sigma_4$   & octahedral  &  $\{  \pm \bu_1 ,  \pm \bu_2 , \pm \bu_3 \}$  &  orthogonal axes      \\
             &  = cubic &   &   \\
 $Y  = A_5$   & icosahedral   &  $\{ \pm \bu_1 ,  \dots , \pm \bu_6 \}$   & 
 $| \bu_i ^\T \bu_j |= 5^{-1/2}$ for $ i \ne j$  \\
            & = dodecahedral   &      &    \\
\hline
       \end{tabular}
       }
\\ ~ \\
{\scriptsize  The $\bu_i$ are unit vectors.}
\end{center}
 \label{table: frame}
 \end{table}

Special cases of Table 1 include the 7 crystal systems: 
triclinic, mono\-clinic,  trigonal, tetragonal, orthorhombic, hexagonal and cubic with sym\-metry groups 
$C_1$, $C_2$, $C_3$,  $C_4$, $D_2$, $D_6$ and $O$, respectively.

\section{Transforming symmetric frames to tensors}

\subsection{Embeddings of the sample spaces}

In order to carry out statistics on $SO(3)/K$, we shall take the embedding 
approach used in, e.g.\  \S10.8 of Mardia \& Jupp (2000). 
We shall embed  $SO(3)/K$ in an inner-product space, $E$, on which $SO(3)$ acts.
 The embedding will be a well-defined equivariant one-to-one function 
$\bt : SO(3)/K \rightarrow E$
such that $\bt ([\bU])$ has expectation $\bzero$ if $[\bU]$ is uniformly distributed on $SO(3)/K$.
For $E = L^2(SO(3))$, the space of square-integrable functions on $SO(3)$, 
a very wide class of such embeddings can be obtained by averag\-ing over $K$.
Let $\bt_0 : SO(3) \to  L^2(SO(3))$ be an embedding used in the Hilbert space approach 
to Sobolev tests of uniformity;  see \S\S 10.8, 13.2.2 of Mardia \& Jupp (2000), Gin\'e (1975)   
and \S4 of Prentice (1978).  
Define $\bt : SO(3)/K \to  L^2(SO(3))$ by $\bt ([\bU]) = |K|^{-1} \sum_{\bR  \in K} \bt_0(\bU \bR)$, 
where $|K|$ denotes the number of elements in $K$. 
If $\bt$ is one-to-one then it is an embedding.
In general, such $\bt$ are quite complicated, so in this paper, for each $K$, we focus on 
a simple choice of embedding, $\bt_K$, of $SO(3)/K$
 into an appropriate space of symmetric tensors.  
These $\bt_K$ are given in Table 2.   
Corresponding expressions for $\langle \bt_K([\bU]),   \bt_K([\bW]) \rangle$ with $\bU, \bW \in SO(3)$ 
are given in Table 3.      
Here $\langle\cdot ,   \cdot \rangle$ is the standard inner product on the appropriate tensor product. 

\begin{table}[ht]
\begin{center}
\caption{Some embeddings $\bt_K :SO(3)/K \rightarrow  E$.} 
{\scriptsize  
\begin{tabular}{llll}
\hline 
	Group, $K$     &   $\bt_K$          \\ 
\hline
   $C_1$       &  
$\bt_{C_1}(\bu_1, \bu_2, \bu_3) =  (\bu_1, \bu_2, \bu_3)$     \\ 
   $C_2$           &  
$\bt_{C_2}(\bu_0, \pm \bu_1) 
=  \left( \bu_0,  \bu_1 \bu_1^\T - (1/3) \bI_3  \right)$    \\  
    $C_r \quad  (r \ge 3)$   &        \\
  $r$ odd  &  
$\bt_{C_r}([ \bu_1 ,  \dots , \bu_r]) = 
\left( \bu_0, \sum_{i=1}^r \otimes ^r \bu_i  \right)$ \\  
  $r$ even  &  
$\bt_{C_r} ([ \bu_1 ,  \dots , \bu_r]) = 
\left( \bu_0, \sum_{i=1}^r  \otimes ^r \bu_i  
 - \, r/(r+1) \,  \mathrm{symm} ( \otimes ^{r/2} \bI_3) \right)$ \\    
 $D_2$  & 
  $\bt_{D_2}( \pm \bu_1, \pm \bu_2) = (\bu_1 \bu_1^\T - (1/3) \bI_3, \bu_2 \bu_2^\T - (1/3) \bI_3 ,\bu_3 \bu_3^\T - (1/3) \bI_3)$   \\  
$D_r \quad  (r \ge 3)$  &  \\
   $r$ odd   &  $\bt_{D_r}([ \bu_1 ,  \dots , \bu_r]) =   \sum_{i=1}^r  \otimes ^r \bu_i   $   \\  
$r$ even  &  
$\bt_{D_r} ([ \bu_1 ,  \dots , \bu_r]) = 
 \sum_{i=1}^r  \otimes ^r \bu_i  
 -  \, r/(r+1) \,  \mathrm{symm} ( \otimes ^{r/2} \bI_3)$ \\  
 $T$    
&  $\bt_{T}( \{  \bu_1 ,  \bu_2 , \bu_3 , \bu_4  \}  ) 
  =  \otimes ^3   \bu_1   +  \otimes ^3  \bu_2  +  \otimes ^3 \bu_3 + \otimes ^3 \bu_4$       \\  
$O$   &  $\bt_O( \{ \pm \bu_1 , \pm \bu_2  , \pm \bu_3 \} )
=  \otimes  ^4 \bu_1   +   \otimes  ^4 \bu_2   +   \otimes  ^4 \bu_3  
 -  (3/5) \, \mathrm{symm} ( \otimes ^2 \bI_3 ) $    \\
  $Y$    & 
$\bt_Y(   \{  \pm \bu_1 ,  \dots , \pm \bu_6  \}  ) =  
 \sum_{i=1}^6 \otimes^{10}  \bu_i - 
\, (6/11)  \mathrm{symm} ( \otimes ^5 \bI_3 )$   \\
\hline
\end{tabular}
}
\\ ~ \\
{\scriptsize  For $C_r$, $\bu_0 = \left\{  \sin (2 \pi / r) \right\}^{-1}  \bu_1  \times \bu_2$.
For $D_2$, $\bu_3 = \pm \bu_1 \times \bu_2$. 
 `symm' denotes symmetrization over permutations of factors of the tensor product.}
 \end{center}
\label{table: embed}
\end{table}

Define $\rho^2$ by
\begin{equation}
 \rho^2 = \| \bt ([\bU]) \|^2 ,
\label{rho.sq}
\end{equation}
which has the same value for all $\bU$ in $SO(3)$.
Then $\bt$ embeds $SO(3)/K$ in the sphere of radius $\rho$ with centre the origin in the vector space $E$.

\begin{table}[ht]
\begin{center}
\caption{Inner products of transforms of symmetric frames. }  
{\scriptsize  
\begin{tabular}{ll}
	Group, $K$     &   Inner product          \\
\hline
   $C_1$    &  $\langle \bt_{C_1}( \bu_1 ,  \bu_2 , \bu_3 ),  \bt_{C_1}( \bv_1 ,  \bv_2 , \bv_3) \rangle
 =  \bu_1^\T \bv_1  + \bu_2^\T \bv_2 + \bu_3^\T \bv_3$   \\
 $C_2$   &  $\langle \bt_{C_2}( \bu_0 ,  \pm \bu_1),  \bt_{C_2}( \bv_0 ,  \pm \bv_1) \rangle 
 = \bu_0^\T \bv_0 + \left( \bu_1^\T \bv_1 \right) ^2 - 1/3$  \\
 $C_r \quad  (r \ge 3)$   &        \\
   $r$ odd   & $\langle \bt_{C_r} ([\bu_1 , \dots,  \bu_r]),  \bt_{C_r} ([ \bv_1 ,  \dots ,  \bv_r])  \rangle
= \bu_0^\T \bv_0 + \sum_{i=1}^r \sum_{j=1}^r  \left( \bu_i^\T \bv_j\right)^r$  \\
     $r$ even   & $\langle \bt_{C_r} ([\bu_1 , \dots,  \bu_r]),  \bt_{C_r} ([ \bv_1 ,  \dots ,  \bv_r])  \rangle
= \bu_0^\T \bv_0 + \sum_{i=1}^r \sum_{j=1}^r  \left( \bu_i^\T \bv_ j\right)^r -  r^2/(r+1)$        \\
   $D_2$  &  $\langle \bt_{D_2}( \pm \bu_1 ,  \pm \bu_2),
  \bt_{D_2}  \pm \bv_1 ,  \pm \bv_2) \rangle
 = \left( \bu_1^\T \bv_1\right) ^2  + \left( \bu_2^\T \bv_2 \right) ^2      
 + \left( \bu_3^\T \bv_3 \right)^2   - 1$   \\
    $D_r \quad  (r \ge 3)$   &       \\
 $r$ odd   &  
$\langle \bt_{D_r} ([\bu_1 , \dots,  \bu_r]), \b t_{D_r} ([ \bv_1 ,  \dots ,  \bv_r])  \rangle  
=  \sum_{i=1}^r \sum_{j=1}^r  \left( \bu_i^\T \bv_j\right)^r$   \\
 $r$ even   & $\langle \bt_{D_r}([\bu_1 , \dots,  \bu_r]),  \bt_{D_r}([\bv_1 ,  \dots ,  \bv_r])  \rangle
=  \sum_{i=1}^r \sum_{j=1}^r  \left( \bu_i^\T \bv_j \right)^r  - r^2/(r+1)$   \\
   $T$   & $\langle  \bt_{T}(  \{  \bu_1 ,  \bu_2 , \bu_3 , \bu_4   \}  ) ,  
   \bt_{T}(  \{  \bv_1 ,  \bv_2 , \bv_3 , \bv_4   \} ) \rangle
  = \sum_{i=1}^4 \sum_{j=1}^4  \left( \bu_i^\T \bv_j \right)^3$   \\
  $O$  &  $\langle \bt_O( \{ \pm \bu_1 ,  \pm \bu_2,  \pm \bu_3 \} ),
 \bt_O ( \{ \pm \bv_1 ,  \pm \bv_2,  \pm \bv_ 3 \} )  \rangle
 = \sum_{i=1}^3 \sum_{j=1}^3  \left( \bu_i^\T \bv_j\right)^4 -  9/5$   \\  
   $Y$   &  $\langle  \bt_{Y}(  \{  \pm \bu_1 ,  \dots , \pm \bu_6  \}  ) , 
    \bt_{Y}(  \{  \pm \bv_1 ,  \dots , \pm \bv_6  \} ) \rangle
  = \sum_{i=1}^6 \sum_{j=1}^6  \left( \bu_i^\T \bv_j\right)^{10}  - 36/11$   \\
\hline
\end{tabular}
}
{\scriptsize  For $C_r$, $\bu_0 = \left\{  \sin (2 \pi / r) \right\}^{-1}  \bu_1  \times \bu_2$.
For $D_2$, $\bu_3 = \pm \bu_1 \times \bu_2$.}
\end{center}
\end{table}

Each symmetric frame can be represented by an element $\bU$ of $SO(3)$.
In the triclinic case, where $K = C_1$, $\bU$ is unique and $\bt ([\bU]) =  \bU$.
We have restricted our attention to point groups, $K$, of the first kind, i.e., excluding reflections.
  However, in situations where reflection symmetries are also present we can
  adopt a right-handed convention for all orientations, and then neglect
reflections.  
For example, we can treat observations on 
$O(3)/\{ \bI_3, - \bI_3 \}$ in the same way as those on $SO(3) = SO(3)/C_1$.

\subsection{Sample mean}

Observations $[\bU_1], \dots , [\bU_n]$ in $SO(3)/K$ can usefully be summarized by the sample mean
$\bar{\bt}$  of their images by $\bt$, i.e., by $\bar{\bt} = n^{-1} \sum_{i=1}^n \bt ([\bU_i])$.
The sample mean $[{\bar \bU}]$ is  defined as the $[\bU]$ in $SO(3)/K$ that maximizes
$\langle \bt( [{\bar \bU}]) , \bar{\bt} \rangle$.
Although $[{\bar \bU}]$ is not necessarily  unique,
it  follows from Theorem 3.2 of Bhatta\-charya \& Patrangenaru (2003)   
that if $[\bU_1], \dots , [\bU_n]$ are generated by 
a continuous distribution then $[{\bar \bU}]$  is unique with probability 1.  

\subsection{Sample dispersion}

A sensible measure of dispersion is 
\begin{equation}
  d = \rho ^2 - \| {\bar \bt} \|^2 , 
\label{d}
\end{equation} 
analogous to the quantity $1 - \bar{R}^2$ used for spherical data; 
see p.\ 164 of Mardia \& Jupp (2000).   
The dispersion satisfies the inequalities $0 \le d \le \rho ^2$, 
where $\rho^2$ is defined in (\ref{rho.sq}).
Since $\bt$ is one-to-one, $d = 0$ if and only if $[\bU_1] = \dots = [\bU_n]$.
Transformation of $[\bU_1], \dots , [\bU_n]$ to 
$[\bV \bU_1], \dots , [\bV \bU_n]$ with $\bV$ in $SO(3)$ leaves $d$ unchanged.
If $K = C_1$ then $d = 3 - \mathrm{trace} \left( {\bar \bR}^2 \right)$,
where ${\bar \bR} = \left( {\bar \bX} ^T {\bar \bX}   \right) ^{1/2}$ with ${\bar \bX} = {\bar \bt}$, 
the sample mean of $\bX = (\bu_1 , \bu_2, \bu_3)$, as in p.\ 290 of Mardia \& Jupp (2000).
If $K = D_2$ then $d = d_1$, where $d_1$ is one of the measures of  
dispersion defined in \S 2.3 of Arnold \& Jupp (2013).

\section{Tests of uniformity}

\subsection{A simple test}

The uniform distribution on $SO(3)/K$ is the unique distribution that is invariant under the action 
of $SO(3)$ on $SO(3)/K$ in which $\bV$ in $SO(3)$ maps $[\bU]$ to $[\bV \bU]$.
Since the embeddings $\bt$ were chosen so that  ${\rm E} \left\{ \bt ([\bU] )\right\} = \bzero$ for $\bU$
uniformly distributed on $SO(3)/K$,  
it is intuitively reasonable to reject uniformity if ${\bar \bt}$ is far from $\bzero$,
i.e.\ if $n \left\| {\bar \bt} \right\| ^2$ is large.
Significance can be assessed using simulation from the uniform distribution.
For large samples, the following asymptotic result can be used.

\medskip
\noindent
{\bf  Proposition 1}

\emph{
Given a random sample on  $SO(3)/K$, define $S$ by 
\begin{equation}
    S = (\nu/\rho^2) n \left\| {\bar \bt} \right\| ^2 
 = n \nu (1 -  d/\rho^2) ,
       \label{S}
\end{equation}
where $\rho^2$ and $d$ are given by (\ref{rho.sq}) and (\ref{d}), respectively, and $\nu$ is the dimension of $E$. 
\begin{enumerate}
\item[(i)] For $K = C_1, D_r$ with $r \ge 2$, $T, O$ or $Y$, 
under uniformity, the asymptotic distribution of $S$ is
$S  \sim   \chi_{\nu}^2$, as $n \rightarrow \infty$.
\item[(ii)] For $K = C_2$, 
\[
  S  = (1/3) S_R + (2/15)S_B ,
\]
where $S_R = 3 n {\bar R}^2$ is the Rayleigh statistic for uniformity of $\bu_0$ and 
$S_B = (15/2)n \left\{ \mathrm{tr} ({\bar \bT} ^2) - (1/3) \right\}$ is the Bingham statistic for uniformity of $\pm \bu_1$, ${\bar R}$ being the mean resultant length of $\bu_0$ and ${\bar \bT}$ being the sample scatter matrix of $\pm \bu_1$.
Under uniformity, $S_R$ and $S_B$ are asymptotically independent with asymptotic distributions 
$\chi_3^2$ and $ \chi_5^2$, respect\-ively.
\item[(iii)] For $K = C_r$ with $r \ge 3$, 
\[
   (\nu_C/\rho_C^2) S  =   (1/3) S_R  + (\nu_D/\rho_D^2) S_D ,
\]
where the subscripts $C$ and $D$ refer respectively to $C_r$-frames and the corresponding $D_r$-frames obtained by replacing 
the directed normal to the plane of a 
$C_r$-frame by the undirected normal, and  $S_R = 3 n {\bar R}^2$ is the Rayleigh statistic for uniformity of $\bu_0$.
Under uniformity, $S_R$ and $S_D$ are asymptotically independent with asymptotic distributions 
$\chi_3^2$ and $ \chi_{\nu_D}^2$, respectively.
\end{enumerate}
}

\medskip

Values of $\rho^2$ and $\nu$ are given in Table 4.   
In the case $K = C_1$, $S$ is  the Rayleigh statistic 
(Mardia \& Jupp, 2000, p.\ 287) for testing uniformity on $SO(3)$. 
In the case $K = D_2$, $S$ is the statistic given in 
\S 3 of Arnold \& Jupp (2013) for testing uniformity on $O(3)/\mathbb{Z}_2^3$.

\begin{table}[ht]
\begin{center}
\caption{Values of squared radius, $\rho^2$, and dimension, $\nu$}  
{\small  
\begin{tabular}{lcr}
	Group     & $\rho^2$   &   $\nu$           \\
\hline
   $C_1$        & 3 &  9  \\
    $C_2$     & 5/3 &  8    \\ 
 $C_r \quad (r \ge 3)$    &    &  \\
  $r$ odd     &  $1 +2^{1-r} r^2$   &   $(r+2)(r+1)/2 + 3$  \\
  $r$ even     &  $1 + r^2 2^{1-r}  \left\{ 1+ 2^{-1}  {r \choose r/2} \right\}    - r^2/(r+1)$  &   $(r+2)(r+1)/2 + 3$    \\
   $D_2$     & $2$ & 10   \\
   $D_r \quad (r \ge 3)$     &    &   \\
  $r$ odd     &   $2^{1-r} r^2 $  &  $(r+2)(r+1)/2$   \\
  $r$ even    &  $r^2 2^{1-r}  \left\{ 1+ 2^{-1}  {r \choose r/2} \right\} - r^2/(r+1)$   & $(r+2)(r+1)/2$   \\
     $T$    & 32/9   &  10 \\
    $O$                & 6/5  &  9   \\
    $Y$                &  $18816/6875$  &  21  \\
\hline
       \end{tabular}
       }
\end{center}
\end{table}

\subsection{Some consistent tests}

The test of uniformity based on $S$ is consistent only against alternatives for 
which ${\rm E} \{ \bt ([\bU]) \}$ is non-zero.
For example, in any equal mixture of two frame cardioid distributions with densities  
(\ref{cardioid}) having concentrations $\kappa$ and $ - \kappa$, $E \{ \bt ([\bU]) \} = \bzero$, and so, 
in asymptotically large samples, $S$ cannot distinguish between such mixtures and the uniform distribution.

Tests of uniformity on $SO(3)/K$ that are consistent against all alternatives can be obtained as 
follows by averaging over $K$ Prentice's generalization to $\mathbb{R}P^3$ 
of Gin\'e's $G_n$ test of 
uniformity; see \S4 of Prentice (1978) and \S13.2.2 of Mardia \& Jupp (2000). 
Given $[\bU_1], \dots , [\bU_n]$ in $SO(3)/K$,  
 with representatives $\bU_1, \dots , \bU_n$ in $SO(3)$, put
\begin{equation}
T_G = -  \sum_{i=1}^n \sum_{j=1}^n  \sum_{\bR \in K}
\left\{ 3 - {\rm trace} \left( \bU_i ^\T \bU_j \bR \right) \right\} ^{1/2}   ,
\label{T_G}
\end{equation}
cf.\ the construction in \S 2 of Jupp \& Spurr (1983).    
Uniformity is rejected if $T_G$ is large compared with the randomization
distribution obtained by replacing $\bU_1, \dots , \bU_n$ by $\bR_1 \bU_1, \dots , \bR_n \bU_n$, 
where $\bR_1, \dots , \bR_n$ are independent random rotat\-ions obtained from 
the uniform distribution on $SO(3)$.
The ortho\-rhombic case, $K = D_2$, is considered in \S 3 of Arnold \& Jupp (2013).
It follows from Theorem 3.1 of Jupp \& Spurr (1983)  
and  the consistency of Gin\'e's test on $\mathbb{R}P^3$ (Mardia \& Jupp, 2000, p.\ 289)
that the test based on $T_G$ is consistent against all alternatives.
More general Sobolev statistics on $SO(3)/K$ can be obtained from Sobolev statistics on 
$SO(3)$ by averaging over $K$, as in (\ref{T_G}).
 
Permutational multi-sample tests, tests of symmetry,  tests of independence, 
and goodness-of-fit 
tests for symmetric frames can be obtained by applying the machinery of 
Wellner (1979), Jupp \& Spurr (1983), Jupp \& Spurr (1985) and Jupp (2005),  
respectively, to the embedding $\bt$.  
These tests of independence are considered in \S 7.

\section{Distributions on $SO(3)/K$}

\subsection{A general class of distributions}

An appealing class of distributions on $SO(3)/K$ consists of those with 
 densities of the form 
\begin{equation}
f([\bU]; [\bM], \kappa ) = g \left( \langle \bt([\bU]) ,   \bt([\bM]) \rangle ; \kappa \right) ,
\label{UARS}
\end{equation}
where $g \left( \cdot ; \kappa \right)$ is a suitable known function and $[\bM] \in SO(3)/K$.
The para\-meter $[\bM]$ measures location and $\kappa$ measures concentration.
If $g \left( \cdot ; \kappa \right)$ is a strictly increasing function,
as in (\ref{Watson}) or Le\'on et al.\ (2006) 
with $\kappa > 0$, then the mode is $[ \bM]$. 

In the case $K = C_1$, the densities (\ref{UARS}) depend on $\bU$ only through 
$\mathrm{trace} ( \bU \bM^\T)$ and the axes and the rotation angles of the random rotations 
are independent, with the axes being uniformly distributed.
These distributions were introduced by Bingham et al.\ (2009) 
under the name of uniform axis-random spin distributions and by Hielscher et al.\ (2010) under the name of radially symmetric distributions.
For $K \ne C_1$, elements of $SO(3)/K$ do not have well-defined axes and, in general,  the distributions on $SO(3)$ with densities ${\tilde f}$ of the form 
${\tilde f}(\bU) = f([\bU]; [\bM], \kappa )$ do not have uniformly distributed axes.

Taking $g(x; \kappa)$ proportional to $e^{\kappa x}$ in (\ref{UARS}) gives the
densities of the form 
\begin{equation}
f([\bU]; [\bM], \kappa ) = c(\kappa)^{-1}  \exp \{ \kappa \langle \bt ([\bU]) ,  \bt ([\bM]) \rangle \} .
\label{Watson}
\end{equation}
For $\kappa > 0$, the mode is $[\bM]$ and
the maximum likelihood estimate of $[\bM]$ is the sample mean. 
The family (\ref{Watson}) is a subfamily of the crystallo\-graphic exponential family introduced by Boogaart (2002, \S 3.2).
For $K = C_1$, (\ref{Watson}) is the density of the matrix Fisher distribution with para\-meter matrix $\kappa \bM$ and 
$c( \kappa) = {_{0}F_1}(3/2, ( \kappa ^2/4) \bI_3)$ (Mardia \& Jupp, 2000, \S 13.2.3).
For $K = D_2$, (\ref{Watson}) is the density of the equal concentration frame Watson 
distributions considered in Arnold \& Jupp (2013, \S 6.1).
Taking $g(x; \kappa)$ proportional to $\left( 1 + x \right) ^{\kappa}$ in (\ref{UARS}) gives the
densities of the form
\begin{equation}
f([\bU]; [\bM], \kappa ) = c(\kappa)^{-1} 
\left\{ 1 + \langle \bt ([\bU]) ,   \bt ([\bM]) \rangle \right\} ^{\kappa} .
\label{Leonetal}
\end{equation}
For $K = C_1$, these densities are those of the de la Vall\'ee Poussin   
distributions introduced by Schaeben (1997),  
 and, under the name of Cayley distributions, by Le\'on et al.\ (2006).

Taking $g(x) = 1 + \kappa x$ with $0 \le \kappa \le \rho^{-2}$ in (\ref{UARS}) gives the
densities  
\begin{equation}
f([\bU]; [\bM], \kappa ) = 1 + \kappa \langle \bt ([\bU]) ,  \bt ([\bM]) \rangle  
\label{cardioid}
\end{equation}
of the frame cardioid distributions, which are analogous to the cardioid distributions on the circle 
(Mardia \& Jupp, 2000, \S 3.5.5).
Useful estimators of $[\bM]$ and $\kappa$ in (\ref{cardioid}) are 
the moment estimators, $[{\hat \bM}]$ and $\hat{\kappa}$, where $[{\hat \bM}]$ is the sample mean defined in \S 3.2 and
 $\hat{\kappa} = (1 - 1/n)^{-1} s^{-2} \langle\bar{\bt} ,   \bt([{\hat \bM}]) \rangle$ with $s^2$ the sample variance of $\langle \bt ([\bU]) , \bt([{\hat \bM}]) \rangle$.

Distributions on $SO(3)/K$ can be identified  with
distributions on $SO(3)$ that are invariant under the action of $K$. One way of generating
such distributions is to average a given distribution on $SO(3)$ over $K$.  
This averaging construction has been used by Walsh et al.\ (2009) in the orthorhombic case, 
Du et al.\ (2016) in the cubic case,
and by Matthies (1982), Gorelova et al.\ (2014) and Niezgoda et al.\ (2016)   
in the general crystallographic case. 
Because the parameters of the distributions (\ref{UARS})  are readily interpretable and 
the distributions (\ref{Watson}), being exponential models, have pleasant inferential properties, we 
find these models more useful than many models obtained by averaging over $K$, especially as the latter can be quite demanding numerically.

\subsection{Concentrated distributions}

A standard coordinate system on $SO(3)$ is given by the inverse of a restrict\-ion of the 
exponential map $\bS \mapsto \sum _{k=0}^{\infty} (k!)^{-1} \bS^{k}$ 
from the space of skew-symmetric $3 \times 3$ matrices to $SO(3)$.
This can be modified to provide coordinate systems on $SO(3)/K$.
Let $[\bM]$ be an element of $SO(3)/K$. 
There are neighbourhoods 
$\mathcal{N}_{[\bM]}$ of $[\bM]$ in $SO(3)/K$ and $\mathcal{V}$ of $\bzero$ in $\mathbb{R}^3$ such that 
each $[\bU]$ in $\mathcal{N}_{[\bM]}$ can be written uniquely as 
$[\bU] = [\bM  \exp \left\{ \bA (\bv) \right\} ]$,
where 
\[
\bA (\bv) =  
\left(
\begin{array}{ccc}
  o   &  - v_3 & v_2 \\
  v_3   & 0         & - v_1 \\
- v_2  &  v_1 & 0 
\end{array} \right) 
\]
with  $\bv = ( v_1,  v_2,  v_3)^\T$ in $\mathcal{V}$.
Define $p_{[\bM]}$ from $\mathcal{N}_{[\bM]}$ to $\mathcal{V}$ by
$p_{[\bM]}([\bU]) = \bv$, where $[\bU] = [\bM  \exp \left\{ \bA (\bv) \right\} ]$. 
Then $p_{[\bM]}$ is a coordinate system on $\mathcal{N}_{[\bM]}$.
Second-order Taylor expansion about $\bzero$ of $[\bU]$ as a function of $\bv$, together with some computer algebra, gives the high-concentration asymptotic distribution of $[\bU]$.

\medskip
\noindent
{\bf  Proposition 2}

\emph{
For $[\bU]$ near $[\bM]$ in $SO(3)/K$ put $[\bU] = [ \bM \exp \left\{ \bA (\bv) \right\} ]$ for $\bv$ near $\bzero$ in $\mathbb{R}^3$.
If $[\bU]$ has density (\ref{Watson}) with $\bt = \bt_K$ as in Table 2  
then the asymptotic distribution of $\kappa ^{1/2} \bv$
as $\kappa \rightarrow \infty$ is normal
with mean $\bzero$ and variance $\bSigma$,
where $\bSigma$ is given in  Table 5.   
If $[\bU]$ has density (\ref{Leonetal}) with $\bt = \bt_K$ then $(\kappa /2)^{1/2} \bv$ has this asymptotic distribution.
}

\medskip

\begin{table}[ht]
\begin{center}
\caption{High-concentration asymptotic variance, $\bSigma$, of $\kappa ^{1/2}$. } 
{\small  
\begin{tabular}{lc}
	Group     &  $\bSigma$           \\
\hline
  $C_1$            &  $(1/2) \bI_3$   \\
 $C_2$            &  $\mathrm{diag}(1/2, 1/4, 1/6)$       \\ 
 $C_r  \quad (r \ge 3)$     &     
$\mathrm{diag} \left[ (1 + r A_r)^{-1} , (1 + r A_r)^{-1} , \{2 r A_r - r(r-1) A_{r-1} \} ^{-1} \right]$       \\  
   $D_2$  &     $(1/4) \bI_3$    \\  
 $D_r  \quad  (r \ge 3)$     &  
$\mathrm{diag} \left[ ( r A_r)^{-1} , ( r A_r)^{-1} , \{r A_r - r(r-1) A_{r-1} \} ^{-1} \right]$  \\ 
 $T$            &  $0.070 \bI_3$  \\
 $O$           &    $(1/8) \bI_3$     \\
 $Y$            &   $0.026 \bI_3$  \\
\hline
\end{tabular}
}

{\small  
For $C_r$ and $D_r$, $v_3$ is the component of $\bv$ normal to the plane of $\bu_1,  \dots , \bu_r$ and
$A_r = \sum _{k=1}^r  \cos( k 2 \pi/r) ^r$.  }
\end{center}
\end{table}

\section{Tests of location}

\subsection{One-sample tests}

Let $[\bM]$ be an element of $SO(3)/K$ which is some
measure of location of a distribution on $SO(3)/K$.
There are various tests of the null hypothesis that $[\bM] =  [\bM_0]$,
where $[\bM_0]$ is a given element of $SO(3)/K$.
The case with $K = D_2$ was considered by Arnold \& Jupp (2013, \S 8).

Permutation tests can be based on the following symmetries of $SO(3)/K$:
For $\bR$ in $K$, define $\rho_{[\bM_0]} (\bR)$ as the transformation that takes  $[\bU]$ to   
$\rho_{[\bM_0]} (\bR) [\bU] = [ \bM_0 \bR \bM_0^\T \bU]$.
Then $\rho_{[\bM_0]}(\bR)$ is well-defined and preserves $[\bM_0]$. 
  
For a sample summarized by the sample mean ${\bar \bt}$ of $\bt$, 
an appealing measure of the squared distance    
between the sample and $[\bM_0]$  is  \linebreak
$\| {\bar \bt} -   \bt([\bM_0])  \|^2$.
It is appropriate to reject the null hypothesis for large values of  $\| {\bar \bt} -   \bt ([\bM_0])  \|^2 $. 
If the distribution of $[\bU]$ is symmetric under $\rho _{[\bM_0]}$ 
then significance can be assessed by comparing the observed value of 
$\| {\bar \bt} -  \bt ([\bM_0]) \|^2$ 
with its random\-ization distribution, which can be obtained by replacing 
$[\bU_1], \dots , [\bU_n]$ by
$\rho_{[\bM_0]} (\bR_1) [\bU_1], \dots , \rho_{[\bM_0]} (\bR_n) [\bU_n]$,
where $\bR_1, \dots , \bR_n$ are independ\-ent and dis\-tribut\-ed uniformly on $K$.

If $[\bU_1] , \dots , [\bU_n]$ is a sample from a concentrated distribution with 
density (\ref{Watson}) and mode $[\bM]$
then it is sensible to test $H_0 : [\bM] = [\bM_0]$ by applying Hotelling's 1-sample $T^2$ test to 
$p_{[\bM_0]}([\bU_1]), \dots , p_{[\bM_0]}([\bU_n])$, where $p_{[\bM_0]}$ is the projection onto the tangent space given in \S 5.2.

\subsection{Two-sample tests}

Suppose that two independent random samples
 $[\bU_1], \dots , [\bU_n]$ and $[\bV_1], \dots , [\bV_m]$ on $SO(3)/K$
are summarized by the sample means ${\bar \bt}_1$ and ${\bar \bt}_2$ of  
 $\bt ([\bU_1]), \dots , \bt ([\bU_n])$ and $\bt ([\bV_1]), \dots , \bt ([\bV_m])$.  
Then the squared distance 
between the two samples can be measured by
 $\| {\bar \bt}_1 -   {\bar \bt}_2  \|^2$.
It is appropriate to reject the null hypothesis that the parent populations are the same if
$\| {\bar \bt}_1 -   {\bar \bt}_2  \|^2$ is large. 
Significance can be assessed by comparing the observed value of 
$\| {\bar \bt}_1 -   {\bar \bt}_2  \|^2 $ with  its randomization distribution, obtained by sampling from the 
potential values corresponding to the 
partitions of the combined sample into samples of sizes $n$ and $m$.

Suppose that $[\bU_1], \dots , [\bU_n]$ and $[\bV_1], \dots , [\bV_m]$ are samples from concentrated distributions
with density (\ref{Watson}) on $SO(3)/K$.
Let $[{\dot \bM}]$ be the maximum likelihood estimate of the mode $[\bM]$ 
under the null hypothesis that the parent populations are the same.
Then the null hypothesis can be tested by applying Hotelling's 2-sample $T^2$ test to 
$p_{[{\dot \bM}]}([\bU_1]), \dots , p_{[{\dot \bM}]}([\bU_n])$ and  
$p_{[{\dot \bM}]}([\bV_1]), \dots , p_{[{\dot \bM}]}([\bV_m])$, where $p_{[{\dot \bM}]}$ is the projection 
onto the tangent space given in \S 5.2.

\section{Independence, regression and misorientation}

\subsection{Independence}

Let $\bt_j : SO(3)/K_j \rightarrow E$, for $j = 1, 2$, be equivariant functions into some 
common inner-product space $E$
such that $\bt_j ([\bU])$ has expectation $\bzero$ if $[\bU]$ is uniformly distributed on $SO(3)/K_j$.
Then association of random variables $[\bU]$ on $SO(3)/K_1$ and $[\bV]$ on $SO(3)/K_2$ can be measured in terms of association of $\bt_1([\bU])$ 
and $\bt_2([\bV])$.

The general approach of Jupp \& Spurr (1985) leads to the following test of independence.
Given pairs $([\bU_1], [\bV_1]),  \dots , ([\bU_n], [\bV_n])$ in $SO(3)/K_1 \times SO(3)/K_2$, 
independence of $\bU$ and $\bV$ is rejected for large values 
of   \linebreak
$\sum_{i = 1}^n \sum_{j = 1}^n   \langle \bt_1 ([\bU_i]) , \bt_1 ([\bU_j]) \rangle   \langle \bt_2 ([\bV_i]) , \bt_2 ([\bV_j])  \rangle$. 
The observed value of this statistic is compared with the randomization distribution given  
$[\bU_1],  \dots , [\bU_n]$, $[\bV_1],  \dots , [\bV_n]$.  
An alternative randomization test rejects independence for large values of the correlation coefficient $r$ defined in (\ref{r.corr}).
For $K_1 = K_2 = C_1$ and $\bt_1 = \bt_2 = \bt_{C_1}$ of Table 2, 
this is one of the tests considered by Rivest \& Chang (2006).   
	
\subsection{Regression}

A reasonable model for homoscedastic regression of $[\bV]$ in $SO(3)/K_2$ on $[\bU]$ in  
$SO(3)/K_1$ has regression function $[\bU] \mapsto [\bA \bU]$ for some $\bA$ in $SO(3)$ and error distribution that is a mild generalization of 
(\ref{Watson}), so that the density of $[\bV]$ given $[\bU]$ is
\begin{equation}
f( [\bV] \mid [\bU] ; \bA, \kappa ) = c(\kappa)^{-1}  \exp \{ \kappa \langle \bt_2 ([\bV]),  \bt_1([\bA \bU]) \rangle \} . 
\label{homosced}
\end{equation} 
For $K_1 = K_2 = C_1$ and $\bt_1 = \bt_2 = \bt_{C_1}$ of Table 2 , 
model (\ref{homosced}) is a generalization of the spherical regression model of Chang (1986).
It is the submodel $\bA_2 = \bI_3$ of the models with regression function $\bU \mapsto \bA \bU \bA_2$ that
were introduced by Prentice (1989) and explored by Chang \& Rivest (2001) and Rivest \& Chang (2006).   
For $K_1 \ne C_1$, it is not possible in general to extend the model (\ref{homosced}) to have 
regression function of the form $[\bU] \mapsto [\bA \bU \bA_2$].
If $K_1 \ne K_2$ then the $\bt_j$ given in Table~\ref{table: embed} 
are not suitable, since in many cases $ \langle \bt_2 ([\bV]),  \bt_1 ([\bA \bU]) \rangle = 0$ for 
all values of $[\bU]$, $[\bV]$ and $\bA$.  
Instead,
it is sensible to take $E = L^2(SO(3))$ and the $\bt_j$  obtained by the averaging construction described in the first paragraph of \S 3.1.

For $\kappa > 0$, the maximum likelihood estimate of $\bA$ is
\begin{equation}
\hat{\bA}=  \arg \max_{\bA \in SO(3)} \sum_{i=1}^n   \langle  \bt_2([\bV_i])  ,  \bt_1( [\bA \bU_i ]) \rangle .
\label{misorn}
\end{equation} 
In general, $\hat{\bA}$ is a well-defined element of $SO(3)$, rather than an element of some quotient. 

Put $\rho_{12} = \max_{\bU \in SO(3)} \langle \bt_1 ([\bU]), \bt_2 ([\bI_3]) \rangle$.
If $\rho_{12} > 0$ then define 
\begin{equation} 
r =  (n \rho_{12})^{-1}   \sum_{i=1}^n   \langle \bt_1 ([{\hat \bA} \bU_i]) ,  \bt_2 ([\bV_i])  \rangle,
\label{r.corr}
\end{equation} 
where ${\hat \bA}$ is given by (\ref{misorn}).
Then $-1 \le \,  r \le 1$ and $r$ can be regarded as a form of uncorrected sample correlation of   
$[\bU]$ and $[\bV]$.
If $K_1 = K_2 = D_2$, $n = 1$ and $\bt_1 = \bt_2 = \bt$ is defined by $\bt ([\bU]) = \bU \, \mathrm{diag} (1,0,-1) \bU^\T$ 
then $r =\cos \omega$, where $\omega$ is the misorientation angle for $D_2$ introduced by Tape \& Tape (2012).                  
Application of Proposition 2 to the decomposition
\begin{eqnarray*} 
  &  &  \sum_{i=1}^n \left\{ \rho^2  - \langle \bt ([\bI_3])  , \bt ([\bV_i ^\T \bA \bU_i]) \rangle \right\}   \\
 & = & \sum_{i=1}^n \left\{ \rho^2  - \langle \bt ([\bI_3])  , \bt ([\bV_i ^\T \hat{\bA} \bU_i]) \rangle \right\}   \\
 &  &  + \sum_{i=1}^n \left\{ \langle \bt ([\bI_3])  , \bt ([\bV_i ^\T \hat{\bA} \bU_i]) \rangle 
     - \langle \bt ([\bI_3])  , \bt( [\bV_i ^\T \bA \bU_i]) \rangle   \right\}   \qquad 
\end{eqnarray*}     
gives the following high-concentration asymptotic distributions.

\medskip
\noindent
{\bf  Proposition 3}

\emph{
For $(\bU_1, \bV_1), \dots , (\bU_n, \bV_n)$ from model (\ref{homosced}),
\begin{enumerate}
\item[(i)]  
Asymptotically, for large $\kappa$,
\begin{eqnarray}
2 \kappa \, \sum_{i=1}^n \left\{ \rho^2  - \langle \bt ([\bI_3])  , \bt ([\bV_i ^\T \bA \bU_i]) \rangle   
 \right\} & \sim & \chi^2_{3n} , \nonumber  \\
2 \kappa \, \sum_{i=1}^n \left\{ \rho^2  - \langle \bt ([\bI_3])  , \bt ([\bV_i ^\T \hat{\bA} \bU_i]) \rangle    
 \right\}  & \sim &  \chi^2_{3(n-1)} ,
\label{chisq3(n-1)}  \qquad \\
2 \kappa \, \sum_{i=1}^n \left\{ \langle \bt ([\bI_3])  , \bt ([\bV_i ^\T \hat{\bA} \bU_i]) \rangle 
     - \langle \bt ( [\bI_3])  , \bt ([\bV_i ^\T \bA \bU_i]) \rangle   
 \right\} & \sim & \chi^2_3  .  \qquad \qquad 
\label{chisq3}
\end{eqnarray}
and the quantities in (\ref{chisq3(n-1)}) and (\ref{chisq3}) are asymptotically independent.
\item[(ii)] An approximate high-concentration $100 (1 - \alpha) \%$ confidence region for $\bA$ is
\[
\left\{ \bA :  2 \hat{\kappa} \, \sum_{i=1}^n \left\{ \langle \bt ([\bI_3])  , \bt ([\bV_i ^\T \hat{\bA} \bU_i]) \rangle 
     - \langle \bt ([\bI_3])  ,\bt ([\bV_i ^\T \bA \bU_i]) \rangle   \
 \right\} < \chi^2_{3; \alpha}  \right\}  .
\]
\end{enumerate}
}

\medskip

\subsection{Misorientation}

The relationship between ambiguous rotations 
$[\bU]$ in $SO(3)/K_1$ and $[\bV]$ in $SO(3)/K_2$ can be described by the mis\-orient\-ation, which is an element of the double coset space 
$K_1 \mbox{\textbackslash} SO(3)/K_2$.
Since $[\bU]$ and $[\bV]$ are images of  $\bU \bR_1$ and $\bV \bR_2$ 
in $SO(3)$ for any $\bR_1$ in $K_1$ and $\bR_2$ in $K_2$, $(\bU \bR_1)^T \bV \bR_2$ determines a well-defined element of 
$K_1 \mbox{\textbackslash} SO(3)/K_2$. See p.\ 274 of Morawiec (1997).  
In crystallography it is usual to identify $K_1 \mbox{\textbackslash} SO(3)/K_2$ with an asymmetric domain, a neighbourhood of $\bzero$ in
$\mathbb{R}^3$ that is in one-to-one correspondence, modulo null sets, with $K_1 \mbox{\textbackslash} SO(3)/K_2$ under 
$\bv \mapsto \exp \{ \bA (\bv) \}$ followed by projection of $SO(3)$ to  $K_1 \mbox{\textbackslash} SO(3)/K_2$.
Then the mis\-orient\-ation between $[\bU]$ and $[\bV]$ is taken as the element, $\bP$, of the asymmetric domain that satisfies 
$\bV = \bU \bP$ and has smallest rotation angle among all such rotations in the domain. 
In the case in which the conditional distribution of $[\bV]$ given $[\bU]$ is uniform the distributions of the angle and axis of  the mis\-orient\-ation are given in Morawiec (2004, Ch.\ 7).  
For general pairs $([\bU_1], [\bV_1]),  \dots , ([\bU_n], [\bV_n])$ in $SO(3)/K_1 \times SO(3)/K_2$,  
we define the mean mis\-orient\-ation as the element $\hat{A}$ of $SO(3)$ defined in (\ref{misorn}).
An alternative definition of the mean mis\-orientation is the element 
$(\hat{\bA}_1, [\hat{\bA}_2])$ of $SO(3) \times (K \mbox{\textbackslash} SO(3))$ that maximizes  
$\sum_i  \max_{\bR_i \in K} \langle  \bt_1 ([ \bA_1 \bU_i]),  \bt_2 ([\bV_i \bR_i \bA_2]) \rangle$.

\section{Example}

To illustrate the estimators and tests introduced above, we
consider some samples of orientations of diopside crystals.
These crystals are monoclinic, so we can represent their orientations by $C_2$-frames.

The stereonet in Fig.\ 1 shows the $\ \bu_0$ vectors and $\pm \bu_1$ axes
given by orientations of 100 diopside crystals.  
A randomization test of uniformity based on $S$ from (\ref{S}) in  \S4.1
has $p-$value less than $0.001$, leading to decisive rejection of uniformity.

\begin{figure}[!ht]
\begin{center}
\includegraphics[width = 0.75 \hsize]{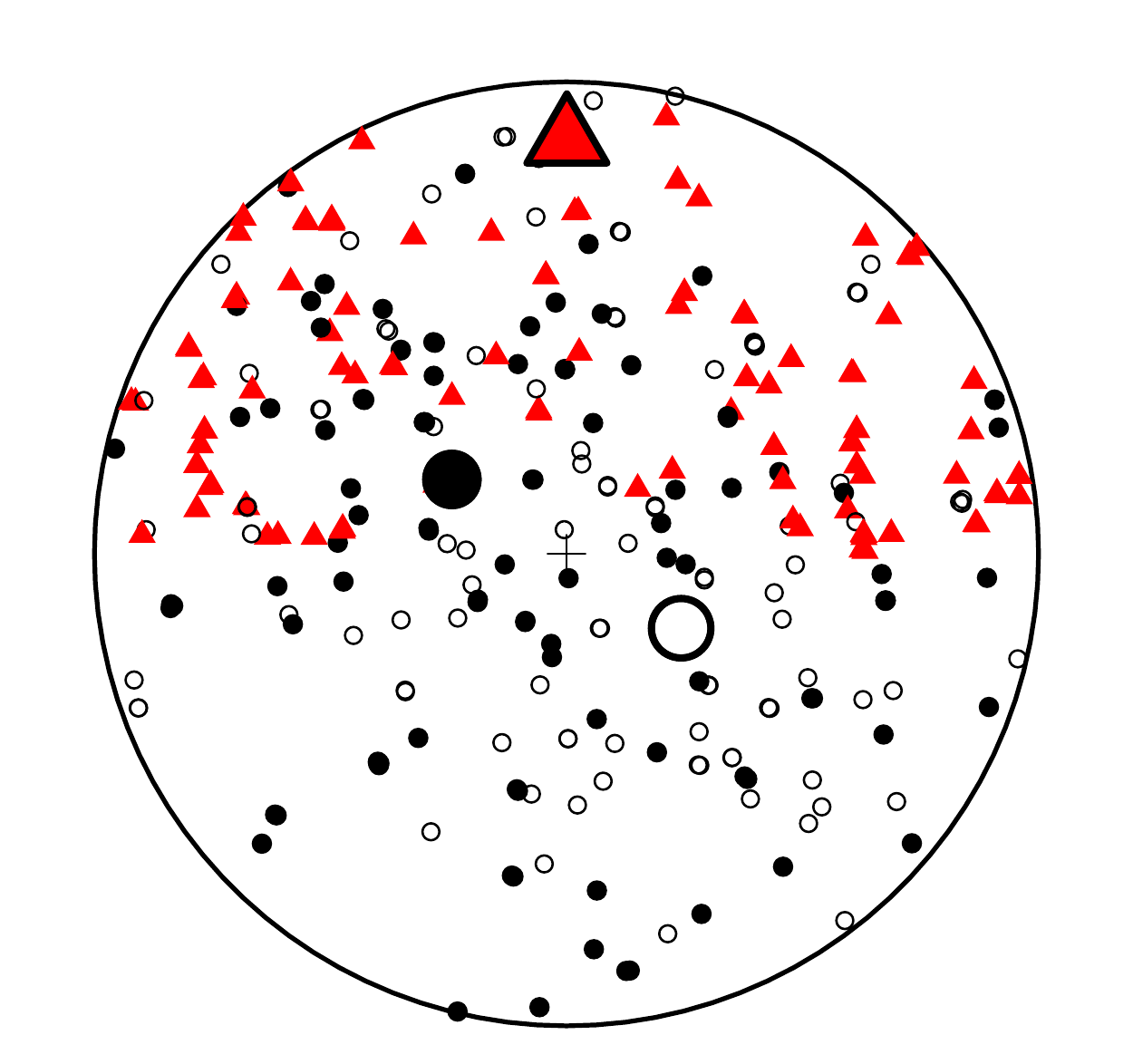}
\end{center}
\caption{Stereonet of $\bu_0$ vectors, shown as red 
triangles, and $\pm \bu_1$ axes, shown as  circles, 
given by orientations of 100 diopside crystals.  
The disc is a stereographic projection of these vectors and axes, showing
the whole of the sphere, so that each axis appears twice with filled circles denoting the 
lower ends of axes and open circles the upper ends.  
The sample mean is shown in large symbols.}
\label{fig:Figure1}
\end{figure}

\begin{figure}
\begin{center}
 \includegraphics[width = 0.75 \hsize]{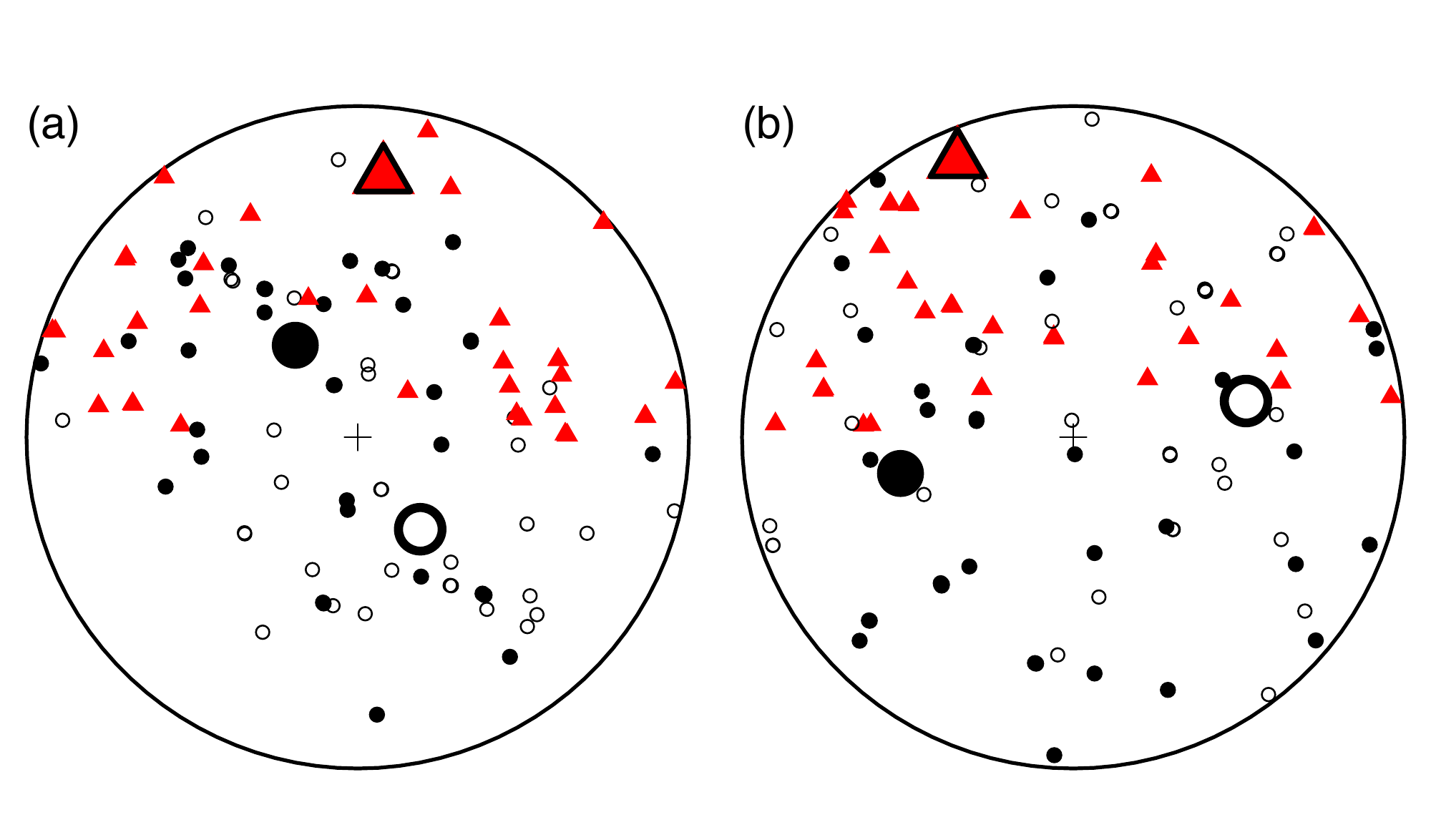}
\end{center}
\caption{Stereonets of $\bu_0$ vectors, shown as red triangles, and $\pm \bu_1$ axes,
shown as circles, given by orientations of two samples of diopside crystals.  
Each disc is a stereographic projection of these axes and vectors, showing
the whole of the sphere, so that each axis appears twice.
(a): 34 orientations from one region of a specimen.
(b): 37 orientations from another region.
The sample means are shown as large symbols.}
\label{fig:Figure2}
\end{figure}

The stereonets in Fig.\ 2 show the  
$\ \bu_0$ vectors and $\pm \bu_1$ axes given by orientations of
34 crystals from one region of a specimen and  37  crystals from another region.  
The two-sample permutation test of \S 6.2 yields a $p$-value of $0.07$ for
equality of the populations of the orientations in the two regions, so the hypothesis
of equality is not rejected.

\section*{Acknowledgements}

We thank David Mainprice for providing the diopside data.

\end{document}